\documentclass[a4paper,12pt,final, reqno]{amsart}
\usepackage{times, a4wide,mathrsfs,amssymb,dsfont, enumerate, xypic}
\usepackage{hyperref}
\usepackage[usenames,dvipsnames]{xcolor}
\hypersetup{colorlinks=true,citecolor=NavyBlue,linkcolor=Brown,urlcolor=Orange}

\usepackage{stmaryrd}
\usepackage{amsmath}
\usepackage{amsfonts}
\usepackage{amsthm}
\usepackage{mathtools}
\usepackage{txfonts}
\usepackage{hyperref}
\usepackage{graphicx,amssymb,amsfonts,amsmath,amscd}
\usepackage{textcomp}
\usepackage[all,ps,cmtip]{xy}
\usepackage{amscd}
\usepackage{xspace}
\usepackage{mathrsfs}
\usepackage{ dsfont }
\usepackage[varg]{pxfonts}

\usepackage{enumerate}
\usepackage{xcolor}
\usepackage{tikz}\usetikzlibrary{decorations.markings,matrix,arrows}
\usepackage{extarrows}

\newcommand{\C}{\mathbb{C}}
\newcommand{\Z}{\mathbb{Z}}

\newcommand{\QQ}{\mathbb{Q}}

\newcommand{\PP}{\mathbb{P}}

\newcommand{\XX}{\mathcal X}

\newcommand{\FF}{\mathcal F}

\newcommand{\CH}{\operatorname{CH}}
\newcommand{\GDCH}{\operatorname{GDCH}}

\newcommand{\YY}{\mathcal Y}
\newcommand{\FFF}{\mathcal F}

\newcommand{\ZZZ}{\mathcal Z}

\newcommand{\pic}{\hbox{Pic}}

\newcommand{\wt}{\widetilde}
\newcommand{\ima}{\operatorname{im}}
\newcommand{\rom}{\romannumeral}

\newcommand{\Gr}{\hbox{Gr}}

\newif\ifHideFoot
\HideFoottrue  
\HideFootfalse  

\ifHideFoot

\newcommand{\Robert}[1]{}
\newcommand{\Charles}[1]{}

\else 

\newcommand{\marg}[1]{\normalsize{{
			\color{red}\footnote{{\color{blue}#1}}}{\marginpar[\vskip
			-.25cm{\color{red}\hfill$\Rightarrow$\tiny\thefootnote}]{\vskip
				-.2cm{\color{red}$\Leftarrow$\tiny\thefootnote}}}}}
\newcommand{\Robert}[1]{\marg{(Robert) #1}}
\newcommand{\Charles}[1]{\marg{(Charles) #1}}

\fi

\newtheorem{theorem}{Theorem}[section]

\newtheorem{lemma}[theorem]{Lemma}
\newtheorem{corollary}[theorem]{Corollary}
\newtheorem{proposition}[theorem]{Proposition}
\newtheorem{conjecture}[theorem]{Conjecture}
\newtheorem{thm2}{Theorem}
\newtheorem{conj2}{Conjecture}
\newtheorem{prop2}{Proposition}

\theoremstyle{definition}

\newtheorem{convention}{Conventions}

\newtheorem{remark}[theorem]{Remark}
\newtheorem{definition}[theorem]{Definition}

\newtheorem{nonumberingt}{Acknowledgements}

\begin{document}

	\author[Robert Laterveer]
	{Robert Laterveer}
	\address{Institut de Recherche Math\'ematique Avanc\'ee,
		CNRS -- Universit\'e 
		de Strasbourg,\
		7 Rue Ren\'e Des\-car\-tes, 67084 Strasbourg CEDEX,
		FRANCE.}
	\email{robert.laterveer@math.unistra.fr}

	\author[Charles Vial]
	{Charles Vial}
	\address{Fakult\"at f\"ur Mathematik, Universit\"at Bielefeld, Germany}
	\email{vial@math.uni-bielefeld.de}

	\dedicatory{To the memory of Jacob Murre}
	
	\date{\today}

	\title[The BVF conjecture and LLSS eightfolds]{The Beauville--Voisin--Franchetta conjecture\\
		 and LLSS eightfolds}	
	
	\begin{abstract} 
The Chow rings of hyper-K\"ahler varieties are conjectured to
have a particularly rich structure. 
In this paper, we formulate a conjecture that combines the Beauville--Voisin conjecture regarding the subring generated by divisors and the Franchetta conjecture regarding generically defined cycles. 
As motivation, we show that this Beauville--Voisin--Franchetta conjecture for a hyper-K\"ahler variety~$X$ follows from a combination of Grothendieck's standard conjectures for a very general deformation of $X$, Murre's conjecture (D) for~$X$ and the Franchetta conjecture for $X^3$. 
As evidence, beyond the case of Fano varieties of lines on smooth cubic fourfolds, 
we show that this conjecture holds for codimension-2 and codimension-8 cycles on Lehn--Lehn--Sorger--van Straten eightfolds. 
Moreover, we establish that the subring of the Chow ring generated by primitive divisors injects into cohomology.
\end{abstract}

  
\thanks{R.L.\ is supported by Agence Nationale de la Recherche under grant ANR-20-CE40-0023.
The research of Ch.V.\ was funded by the Deutsche Forschungsgemeinschaft (DFG, German Research Foundation) -- SFB-TRR 358/1 2023 -- 491392403}

	\maketitle

	\section*{Introduction}

	Starting with the seminal work of Beauville and Voisin on the Chow ring of K3 surfaces~\cite{BV}, it has been observed that the Chow rings (and more generally the Chow
	motives, considered as algebra objects) of hyper-K\"ahler varieties possess a surprisingly rich structure.
	Our aim is to set the famous Beauville--Voisin conjecture into a wider conjectural framework. 
	Doing so, we formulate a generalized version of the Beauville--Voisin conjecture concerning the Chow ring of hyper-K\"ahler
	varieties, and establish it in the special case of zero-cycles on Lehn--Lehn--Sorger--van Straten eightfolds (LLSS eightfolds, for short).
	These varieties form a 20-dimensional
	locally complete family of hyper-K\"ahler varieties, deformation equivalent to
	the Hilbert scheme of length-4 subschemes on a K3 surface.

	\subsection{The Beauville--Voisin conjecture}
	The following conjecture is by now classical:
	
	\begin{conj2}
		Let $X$ be a hyper-K\"ahler variety. 
		 Consider the $\QQ$-subalgebras
		\[ 	B^\ast(X)\coloneqq \langle \CH^1(X)\rangle \quad \mathrm{and} \quad
		BV^\ast(X)\coloneqq \langle \CH^1(X),  c_j(X)\rangle\]
           of the Chow ring $\CH^*(X)$ with $\QQ$-coefficients
		generated on the one hand by divisors and on the other hand by divisors and
		Chern classes.
		\begin{itemize}
\item(Beauville conjecture) 
The restriction of the cycle class map $B^i(X)\to H^{2i}(X,\QQ)$ is injective
for all $i$.
\item (Beauville--Voisin conjecture) 
The restriction of the cycle class map $BV^i(X)\to H^{2i}(X,\QQ)$ is injective
for all $i$.
		\end{itemize}
	\end{conj2}
	
	The Beauville--Voisin conjecture was first proven (without being stated as such) in the case of K3
	surfaces in the seminal work of Beauville and Voisin \cite{BV}.
	The Beauville conjecture was then formulated by Beauville \cite{Beau3}
	 as an explicit workable consequence of a deeper conjecture stipulating
	the splitting of the conjectural Bloch--Beilinson filtration on the Chow rings
	of hyper-K\"ahler varieties. The Beauville--Voisin conjecture was then stated
	 by Voisin~\cite{V17}
	who established it in the case of Hilbert schemes of points on K3 surfaces of
	low dimension and in the case of Fano varieties of lines on smooth cubic
	fourfolds. The Beauville conjecture has by now been established in many
	cases, including hyper-K\"ahler
	varieties carrying a rational Lagrangian fibration \cite{Riess} as well as Hilbert schemes of points on K3 surfaces \cite{MN}. The Beauville--Voisin conjecture has been established for generalized Kummer
	varieties \cite{FuKummer}, and
	 in addition to the cases mentioned above of K3 surfaces \cite{BV} and Fano varieties of
	lines on cubic fourfolds \cite{V17}, there is a further locally complete family of
	hyper-K\"ahler varieties all of whose members are known to satisfy the
	Beauville--Voisin conjecture, namely that of double EPW sextics; see \cite{Fe, LV} for partial results and \cite{Lat} for a full proof.
	In this note, we give further
	 evidence for the Beauville--Voisin conjecture 
	 by establishing it for 0-cycles on the locally complete family of LLSS eightfolds.
	
	\begin{thm2}\label{main2}
 The Beauville--Voisin conjecture holds for 0-cycles on LLSS eightfolds $Z$. More precisely, there exists a closed point $o_Z \in Z$ such that
 $$BV^8(Z) = \QQ [o_Z] .$$
	\end{thm2}
	
	In fact, we will establish a stronger result; see Theorem~\ref{main3} below.
  As more evidence towards the Beauville conjecture for hyper-K\"ahler varieties, we have: 
	
	\begin{prop2}\label{Beauville-weak2}
		Let $Z$ be an LLSS eightfold. The $\QQ$-subalgebra
		\[  \langle \CH^1(Z)_{\rm prim} \rangle\ \ \subset\ \CH^\ast(Z) \]
		generated by the primitive divisors injects into cohomology.
	\end{prop2}
	
	In view of \cite[Theorem~1.4 and Theorem~4.3(4)]{Voi-triangle}, where Proposition~\ref{Beauville-weak2} is established conditional on its validity in top codimension, 	Proposition~\ref{Beauville-weak2} is a corollary of Theorem~\ref{main2}.
	However, we refer to Proposition~\ref{Beauville-weak} for a more direct proof based on equivariant resolution of indeterminacy and on the existence of  an anti-symplectic involution on $Z$ compatible with Voisin's rational map~\cite[Proposition~4.8]{V14}. 
	
	\subsection{The Beauville--Voisin--Franchetta conjecture}
	While the Beauville conjecture is a direct consequence of the conjectural splitting of the
	Bloch--Beilinson filtration on the Chow ring of hyper-K\"ahler varieties, 
	the Beauville--Voisin conjecture suggests that the Chern classes of hyper-K\"ahler varieties lie in the grade-0 part of this conjectural splitting. 
	More generally, we expect generically defined cycles on hyper-K\"ahler varieties, as defined below, to lie in the grade-0 part of this conjectural splitting. This leads to the Beauville--Voisin--Franchetta Conjecture~\ref{conj2} below.\medskip

	Let $\mathcal F$ be the moduli stack of a locally complete family of polarized hyper-K\"ahler varieties  and let $X$ be any member of the universal family $\mathcal X \to \mathcal F$. 
	For positive integers~$n$, we define the subring of \emph{generically defined cycles} on $X^n$ to be 
	$$\GDCH^\ast_{\mathcal F}(X^n) \coloneqq \operatorname{im} \big( \CH^*(\mathcal X^{\times n}_{/ \mathcal F}) \to \CH^*(X^n) \big).$$
	The generalized Franchetta conjecture for powers of hyper-K\"ahler varieties stipulates that the restriction of the cycle class map 
	$$\GDCH^i_{\mathcal F} (X^n) \to H^{2i}(X^n,\QQ)$$
	is injective for all $i$ and $n$. 
	This conjecture was first stated by O'Grady \cite{OG} for K3 surfaces in the case $n=1$. 
	It was settled for $n=1$ in the case of K3 surfaces of low genus in \cite{PSY, FL}, in the case of Fano varieties of lines on cubic fourfolds in \cite[Theorem~1.9]{FLV} and in the case of LLSS eightfolds first in \cite[Theorem~1.11]{FLV} for 0-cycles and for codimension-2 cycles and then in general in \cite{FLV2}. It was also settled for low values of $n$ and K3 surfaces of low genus in \cite[Theorem~1.5]{FLV}, and for $n=2$ and the Fano variety of lines on a cubic fourfold in \cite[Theorem~1.10]{FLV}.
	It is probably optimistic to expect the generalized Franchetta conjecture to hold for large values of $n$, but it would be interesting in any case to find lower bounds.
	\medskip

As explained in Section~\ref{s:BVF} (see Proposition~\ref{prop:F-MCK} and Corollary~\ref{cor:F-BVF}), a combination of the standard conjectures for the very general fiber of $\mathcal X\to \mathcal F$, of the generalized Franchetta conjecture for powers of $X$, and of Murre's conjecture (D) \cite{Murre} for the generic Chow--K\"unneth decomposition implies the following generalized version of the Beauville--Voisin conjecture:

	\begin{conj2}[Beauville--Voisin--Franchetta conjecture] \label{conj2}
		Let $\mathcal F$ be the moduli stack of a locally complete family of polarized hyper-K\"ahler varieties and let $X$ be any member of the universal family $\mathcal X \to \mathcal F$.
   Consider, for a positive integer $n$, the $\QQ$-subalgebra
	\[ R^\ast(X^n)\coloneqq \big\langle \CH^1(X^n), \  \GDCH^\ast_{\mathcal F}(X^n)\big\rangle\]
	of the Chow ring $\CH^*(X^n)$ with $\QQ$-coefficients
	generated by divisors and
	generically defined cycles.
	Then the restriction of the cycle class map 
	$R^i(X^n)\to H^{2i}(X^n,\QQ)$ is injective
	for all $i$ and all $n$.
\end{conj2}

%

 For $n=1$,	Conjecture~\ref{conj2} holds  for K3 surfaces of low genus by \cite{PSY, FL} and for Fano varieties of lines on smooth cubic fourfolds by combining \cite[Theorem~1.4$(2)$]{V14} and the results of \cite[Section~3]{FLV}. Conjecture~\ref{conj2} also holds for low values of $n$ and for K3 surfaces of low genus by combining \cite{BV} and the results of \cite[Section~5]{FLV}. 
 In this paper, we show that Conjecture~\ref{conj2} holds for 0-cycles and codimension-2 cycles on LLSS eightfolds:
 
 \begin{thm2}\label{main3}
The Beauville--Voisin--Franchetta conjecture holds for 0-cycles and for codimension-2 cycles on LLSS eightfolds $Z$; that is, the cycle class map induces injections
     \[ R^i(Z)\ \hookrightarrow\ H^{2i}(Z,\QQ)\ \ \ \hbox{for}\ i=2\ \hbox{and\ for}\ i=8\ .\] 
 Moreover, there exists a closed point $o_Z \in Z$ such that
$$R^8(Z) = \QQ [o_Z]\ .$$
 \end{thm2}
		
Theorem~\ref{main3}  implies Theorem~\ref{main2}.
Its proof, given in  Section~\ref{S:main}, reduces to the case of the square of the Fano variety of lines $F$  via a result of Chen~\cite{Chen} showing that Voisin's rational map $\psi\colon F\times F \dashrightarrow Z$ can be resolved by blowing up the locus of intersecting lines. 
	The construction of LLSS eightfolds, the Voisin map and the result of Chen are reviewed in Section~\ref{S:LLSS}.
	Since the locus of intersecting lines is in general singular, we have to make use of Fulton's operational Chow cohomology theory as a replacement of the usual Chow groups. Basic facts regarding operational Chow cohomology are gathered in Section~\ref{S:CH}.
	Regarding the case of codimension-2 cycles in Theorem~\ref{main3}, we have to consider Abel--Jacobi trivial cycles on (not necessarily smooth) projective varieties; their definition and basic properties are given in Section~\ref{S:AJ}.\medskip
	
	With this paper, we wish to contribute to the memory of Jacob Murre. We both have fond memories of Murre's person; his gentleness, his modesty and his rigorous precision will continue to be an inspiration.

	\begin{convention} 
		In this note, the word {\sl variety\/} will refer to an integral scheme of finite type over $\C$. 
		For a variety $X$,
		$\CH_j(X)$ will denote the Chow group of $j$-dimensional algebraic cycles on
		$X$ with $\QQ$-coefficients, while $\CH^j(X)$ will denote Fulton--MacPherson's $j$-th operational Chow cohomology group with $\QQ$-coefficients of $X$. For $X$ smooth, $\CH^j(X)$ identifies canonically with the Chow group of $j$-codimensional algebraic cycles on
		$X$ with $\QQ$-coefficients, and 
		for $X$ smooth of dimension $n$ the notations $\CH^j(X)$ and $\CH_{n-j}(X)$ will
		be used interchangeably. 
	\end{convention}

	\begin{nonumberingt} R.L. thanks Ch.V. for an invitation to Bielefeld,
		where this work originated. Thanks to Qingyuan Jiang for kindly answering
		our questions 
		regarding
		the proof of Lemma~\ref{jiang}.
	\end{nonumberingt}

\section{The Beauville--Voisin--Franchetta conjecture}\label{s:BVF}

In this section we explain how, for a hyper-K\"ahler variety $X$, a combination of Murre's conjectures and of the generalized Franchetta conjecture for powers of $X$ implies the Beauville--Voisin--Franchetta Conjecture~\ref{conj2}.

\subsection{Murre's conjectures}
Murre \cite{Murre} introduced the notion of Chow--K\"unneth decomposition
and conjectured that the induced filtration on Chow groups is the conjectural Bloch--Beilinson filtration. 
While the Bloch--Beilinson conjecture involves all smooth projective varieties at once, Murre's conjecture has the advantage that it is stated for each smooth projective variety individually.

\begin{conjecture}[Murre's Conjecture (A)] \label{conj:A}
Any $d$-dimensional smooth projective variety $X$ admits a \emph{Chow--K\"unneth decomposition}, i.e., there exist correspondences $\pi_X^i \in \CH^{d}(X\times X)$, $0\leq i \leq 2d$, such that
\begin{itemize}
\item $\pi_X^i\circ \pi_X^i = \pi_X^i$ for all $i$, and $\pi_X^i\circ \pi_X^j=0$ for all $i\neq j$.
\item $\Delta_X = \pi_X^0 + \cdots + \pi_X^{2\dim X}$ in $\CH^{d}(X\times X)$, where $\Delta_X$ denotes the class of the diagonal.
\item $\pi_X^i$ acts as the identity on $H^i(X,\QQ)$ and as zero on $H^j(X,\QQ)$ for all $i\neq j$.
\end{itemize}
\end{conjecture}

In other words, the Chow motive $h(X)$ of $X$ decomposes as $h(X) = h^0(X) \oplus \cdots \oplus h^{2d}(X)$, where $h^i(X) \coloneqq (X,\pi_X^i)$.
By the K\"unneth formula, Murre's Conjecture (A) is stable under taking products.
 Such a decomposition provides a descending filtration on the Chow groups of $X$ as follows:
$$F^k\CH^i(X) \coloneqq \CH^i\big(h^{2i-k}(X) \oplus \cdots \oplus h^0(X)\big) = \big(\pi_X^{2i-k} + \cdots + \pi_X^0\big)_\ast \CH^i(X).$$

\begin{conjecture}[Murre's Conjectures  (B), (C) and (D)] Let $X$ be a smooth projective variety admitting a Chow--K\"unneth decomposition as in Conjecture~\ref{conj:A} and let $F^\bullet$ be the induced filtration on its Chow groups.
	\begin{enumerate}
		\item[(B)] $F^0\CH^i(X) = \CH^i(X)$ and $F^{i+1}\CH^i(X) = 0$.
		\item[(C)] The filtration $F^\bullet\CH^i(X)$ does not depend on the choice of a Chow--K\"unneth decomposition.
		\item[(D)] $F^1\CH^i(X) = \ker \big(\CH^i(X) \to H^{2i}(X,\QQ) \big)$.
	\end{enumerate}
\end{conjecture}

Murre's Conjectures (A), (B), (C) and (D) hold for smooth projective curves,
and for smooth projective surfaces
Murre~\cite{Murre-surface} established the existence of a Chow--K\"unneth decomposition that satisfies (B) and (D). 
The existence of a Chow--K\"unneth decomposition for abelian varieties was established by Deninger--Murre in \cite{DM}.
 Based on work by Xu--Xu~\cite{XX}, it was established in \cite[Theorem~3]{vial-ab} that Murre's Conjectures (B), (C) and (D) for abelian varieties reduce to Murre's Conjecture (D) for products of curves. 
 More generally, varieties that are Kimura--O'Sullivan finite-dimensional and satisfy Grothendieck's K\"unneth standard conjecture admit a Chow--K\"unneth decomposition by \cite[Proposition 5.3]{J2} and \cite[Proposition~7.5]{Kim1}.
 Further examples of smooth projective varieties that are known to admit a Chow--K\"unneth decomposition include complete intersections, uniruled threefolds \cite{dAMS}, elliptic modular threefolds \cite{GM99}, threefolds and fourfolds with a nef tangent bundle \cite{Iyer, vial-ab},
   threefolds and fourfolds admitting an MRC quotient of dimension $\leq 2$ \cite{Via15}, 
 Fano varieties of lines on smooth cubic fourfolds \cite{SV, FLV} and more generally on smooth cubic hypersurfaces \cite{FLV2},
 Kuga-Sato varieties of Hilbert modular varieties
 \cite{GHM},
 conic fibrations over a surface \cite{NS} and more generally quadric fibrations over a surface \cite{Via-ACF, Bou},
 varieties with finite-dimensional Chow groups \cite{J2, Kim2} and more generally
 varieties with no non-zero Abel--Jacobi trivial cycle classes \cite{Via-AJ}, certain moduli spaces \cite{dBR}, Hilbert schemes of points on surfaces \cite{dCM}, and generalized Kummer varieties \cite{dCM2, FTV, Xu}.
\medskip

The \emph{Bloch--Beilinson conjecture} stipulates the existence of a descending filtration $F^\bullet$ on the Chow groups of smooth projective varieties with the following properties; 
we refer to \cite[Section~2]{J2} for more details.
\begin{enumerate}
	\item $F^0\CH^i(X) = \CH^i(X)$ and $F^1\CH^i(X) = \ker \big(\CH^i(X) \to H^{2i}(X,\QQ) \big)$.
	\item $F^{i+1} \CH^i(X) = 0$.
	\item $F^\bullet$ is a ring filtration on $\CH^\ast(X)$, i.e., $F^k\CH^i(X) \cdot F^l\CH^j(X) \subseteq F^{k+l}\CH^{i+j}(X)$.
	\item $F^\bullet$ is respected by pullbacks and pushforwards along morphisms.
	\item Assuming Grothendieck's K\"unneth standard conjecture, the graded pieces \linebreak $\mathrm{Gr}_F^k \CH^i(X)$ only depend on the homological motive $h^{2i-k}_{hom}(X)$.
\end{enumerate}
We note that the combination of (3) and (4) is equivalent to requiring the filtration to be compatible with the action of correspondences (since the intersection product is induced by pulling back along the diagonal embedding). Regarding (5), we note that (1) together with (3) and (4) imply that the action of a correspondence on $\mathrm{Gr}_F^k \CH^i(X)$ factors through homological equivalence; (5) then means that the action further factors through $\mathrm{End}\big(h^{2i-k}_{hom}(X)\big)$.
\medskip

The link between Murre's conjectures and the Bloch--Beilinson conjecture is provided by the following:

\begin{theorem}[Jannsen {\cite[Theorem~5.2]{J2}}]
	Murre's conjecture is equivalent to the Bloch--Beilinson conjecture. Precisely:
	\begin{itemize}
		\item 
		If Murre's conjectures (A), (B), (C) and (D) hold for all smooth projective varieties, then the filtration induced by Chow--K\"unneth decompositions is the Bloch--Beilinson filtration. 
		\item 
		Conversely, assuming Grothendieck's K\"unneth standard conjecture, the existence of a Bloch--Beilinson filtration implies that Murre's conjectures (A), (B), (C) and (D) hold for all smooth projective varieties.
	\end{itemize}
 In particular, the Bloch--Beilinson filtration, if it exists, is unique.
\end{theorem}

\subsection{The Franchetta property and the BVF conjecture} 
We aim to show (in Corollary~\ref{cor:F-BVF} below) that a combination of Grothendieck's standard conjectures, of Murre's conjecture (D), and of the generalized Franchetta conjecture for hyper-K\"ahler varieties implies the Beauville--Voisin--Franchetta Conjecture~\ref{conj2}.
We start with a definition:

\begin{definition}\label{def:gdc}
	Let $\XX\to B$ be a family of algebraic varieties, where $B$ is smooth, and let $X=X_b$ be the fiber over $b\in B$. We write
	\[ \GDCH^i_B(X)\coloneqq \operatorname{im}\bigl(   \CH^i(\XX)\to \CH^i(X)\bigr) \]
	for the subgroup of {\em generically defined cycles}.
	We say that the family $\XX\to B$ has the \emph{Franchetta property} in degree $i$ if 
	 the restriction of the cycle class map $\GDCH^{i}_B(X) \to H^{2i}(X,\QQ)$ is injective. We say that it has the Franchetta property if it has the Franchetta property in degree $i$ for all $i$.
\end{definition}

By projecting on the first $m$ factors, we observe that for $m\leq n$ the Franchetta property for $\mathcal X_{/B}^n \to B$ in degree $i$ implies the Franchetta property for $\mathcal X_{/B}^m \to B$ in degree $i - (n-m)\dim_B \mathcal X$.
\medskip

The Franchetta property cannot be expected to hold in general. Obviously, if the family $\mathcal X \to B$ is constant (e.g., if $B$ is reduced to a point) and if the cycle class map $\CH^*(X) \to H^*(X,\QQ)$ is not injective (e.g., if $X$ is a smooth projective curve of positive genus), then it does not satisfy the Franchetta property. 
More interestingly,  the universal curve $\mathcal C \to \mathcal M_g$ of genus $g\geq 2$ satisfies the Franchetta property \cite{Har}, 
whereas if $g\geq 4$ its relative square  does not satisfy the Franchetta property: indeed, it is shown in \cite{GG} (cf. also \cite{Yin}) that the generically defined cycle 
  \[ K_C \times K_C - \deg (K_C) \delta_* K_C\ \  \in \CH^2(C\times C)\] 
  is non-trivial, 
where $K_C$ is the canonical divisor of~$C$ and where $\delta\colon C \to C\times C$ is the diagonal embedding.
Nonetheless, for hyper-K\"ahler varieties, with the evidence listed in the introduction, it remains possible that

\begin{conjecture}[Generalized Franchetta conjecture for hyper-K\"ahler varieties \cite{OG, FLV, FLV2}] \label{conj:F}
Let $\mathcal F$ be the moduli stack of a locally complete family of polarized hyper-K\"ahler varieties and let $\mathcal X \to \mathcal F$ be the universal family. Then, for all $n\geq 1$,  	
the family $\mathcal X^n_{/\mathcal F} \to \mathcal F$ has the Franchetta property.
\end{conjecture}

When combined with general conjectures related to algebraic cycles, the Franchetta property has remarkable consequences as we outline below. 
Let $k$ be a field and $\ell$ be a prime invertible in $k$; 
we say that a $d$-dimensional smooth projective variety $X$ over a field~$k$ satisfies \emph{Grothendieck's K\"unneth standard conjecture} if there exist cycles $p^i_X \in \CH^{d}(X\times_k X)$ for $0\leq i \leq 2d$ such that $p^i_X$ acts as the identity on the $\ell$-adic cohomology group $H^i(X_{\bar k},\QQ_{\ell})$ and as zero on $H^j(X_{\bar k},\QQ_{\ell})$ for $j\neq i$. 
If $k=\C$, by Artin's comparison isomorphism, such cycles  $p^i_X$ act as the identity on  $H^i(X,\QQ)$ and as zero on $H^j(X,\QQ)$ for $j\neq i$. 
(The reason for invoking $\ell$-adic cohomology here, and not sticking merely to singular cohomology of complex varieties, will become apparent in the proof of Proposition~\ref{prop:F-CK} below.)
In relation to the existence of a Chow--K\"unneth decomposition, we start with the following easy observation:

\begin{proposition}\label{prop:F-CK}
	Let $\XX\to B$ be a family of $d$-dimensional algebraic varieties, and let $X=X_b$ be the fiber over a closed point $b\in B$.
	Assume the following:
	\begin{enumerate}[(a)]
		\item The very general fiber of $\XX\to B$ satisfies Grothendieck's K\"unneth standard conjecture.
		\item The family $\XX^2_{/B} \to B$ has the Franchetta property in degree $d$.
	\end{enumerate}
	Then $X$ has a unique generically defined Chow--K\"unneth decomposition $\{\pi_X^i \ \vert \ 0 \leq i \leq 2d\}$.
\end{proposition}
\begin{proof}
	First we explain how (a) implies that the K\"unneth projectors for $X$ are algebraic and generically defined. 
	By \cite[Lemma~2.1]{Via-ACF}, the very general fiber of $\XX\to B$ is isomorphic, as abstract schemes, to its geometric generic fiber. Therefore, the algebraic K\"unneth projectors for the very general fiber of $\XX\to B$ provide algebraic K\"unneth projectors
	 for its geometric generic fiber. 
	As in the proof of \cite[Proposition~1.7]{vial-ab}, it follows from the fact that K\"unneth projectors are Galois invariant, that those descend to a set of algebraic K\"unneth projectors
	 for the generic fiber of $\XX\to B$. 
	By spreading, these provide cycle classes in $\CH^d(\XX \times_B \XX)$ 
  and they specialize to algebraic K\"unneth projectors.
Second,	the generalized Franchetta conjecture for $X^2$ then shows that these provide a Chow--K\"unneth decomposition $\{\pi^i_X \ \vert \ 0\leq i\leq 2d\}$ for $X$.
\end{proof}

We note that assumption (a) holds for hyper-K\"ahler varieties that are deformation-equivalent to the Hilbert scheme of length-$n$ subschemes on a K3 surface by \cite{CM}, 
and for  hyper-K\"ahler varieties that are deformation-equivalent to a generalized Kummer variety of dimension $2n$ with $n+1$ prime by \cite{Foster}.
\medskip

Recall from \cite[Section~8]{SV} that a Chow--K\"unneth decomposition $\{\pi_X^i \ \vert \ 0 \leq i \leq 2d\}$ is \emph{multiplicative} if 
  \[ \pi_X^k \circ \delta_X \circ (\pi_X^i \otimes \pi_X^j) =0\ \ \forall k \not= i+j\ ,\]
  where $\delta_X\in\CH^{2d}(X\times X\times X)$ denotes the small diagonal. 
Equivalently, a Chow--K\"unneth decomposition $h(X) = h^0(X) \oplus \cdots \oplus h^{2d}(X)$ is multiplicative if the restriction to $h^i(X) \otimes h^j(X)$ of the intersection product $h(X) \otimes h(X) \to h(X)$ factors through $h^{i+j}(X)$. 
While Chow--K\"unneth decompositions are expected to exist for any smooth projective variety, the existence of a multiplicative  Chow--K\"unneth decomposition is very restrictive. 
For instance, for smooth projective curves, the existence of a multiplicative  Chow--K\"unneth decomposition is equivalent to the vanishing of the Ceresa cycle \cite[Theorem 1.5.5]{Zhang}.
On a positive note, for abelian varieties, the Chow--K\"unneth decomposition of Deninger and Murre \cite{DM} is multiplicative, and, for hyper-K\"ahler varieties, Shen and Vial~\cite[Conjecture~4]{SV} have conjectured based on Beauville's splitting principle~\cite{Beau3} the existence of a multiplicative Chow--K\"unneth decomposition. 
We refer to \cite{FLV-mult} for more details and references on multiplicative Chow--K\"unneth decompositions.
In relation to the latter and in view of Conjecture~\ref{conj:F}, we have the following proposition:

\begin{proposition}\label{prop:F-MCK}
	Let $\XX\to B$ be a family of $d$-dimensional algebraic varieties, and let $X=X_b$ be the fiber over a closed point $b\in B$.
	Assume the following:
	\begin{enumerate}[(a)]
		\item The very general fiber of $\XX\to B$ satisfies Grothendieck's K\"unneth standard conjecture.
		\item The family $\XX^3_{/B} \to B$ has the Franchetta property in degree $2d$.
	\end{enumerate}
Then $X$ has a unique generically defined Chow--K\"unneth decomposition $\{\pi_X^i \ \vert \ 0 \leq i \leq 2d\}$ and it is multiplicative.
\end{proposition}
\begin{proof}
	Since $(b)$ implies that the  family $\XX^2_{/B} \to B$ has the Franchetta property in degree $d$, Proposition~\ref{prop:F-CK}
	provides a generically defined Chow--K\"unneth decomposition $\{\pi^i_X \ \vert \ 0\leq 0\leq 2d\}$ for $X$. We conclude using $(b)$ that it is multiplicative.
\end{proof}

\begin{proposition}\label{prop:F-0}
	Let $\XX\to B$ be a family of $d$-dimensional algebraic varieties, 
	and let $X=X_b$ be the fiber over a closed point $b\in B$.
	Let $n$ be a positive integer.
	Assume the following:
	\begin{enumerate}[(a)]
		\item The very general fiber of $\XX\to B$ satisfies Grothendieck's K\"unneth standard conjecture.
		\item 
		The family $\XX^n_{/B} \to B$ has the Franchetta property.
		\item[(b')]
		 $X$ has a generically defined multiplicative Chow--K\"unneth decomposition $\{\pi_X^i \ \vert \ 0 \leq i \leq 2d\}$.
	\end{enumerate}
	Then, setting
	\begin{equation}\label{eq:CK}
\pi^j_{X^n} \coloneqq \sum_{i_1 + \cdots + i_n = j}\pi_X^{i_1}\otimes \cdots \otimes \pi_X^{i_n}
\quad \mathrm{and} \quad 
\CH^i(X^n)_{(\nu)} \coloneqq \big(\pi^{2i-\nu}_{X^n} \big)_*\CH^i(X^n)\ ,
	\end{equation}
	we have an inclusion of graded $\QQ$-subalgebras
	$$\GDCH^*_B(X^n) \subseteq \CH^*(X^n)_{(0)} \subseteq \CH^*(X^n).$$
\end{proposition}
\begin{proof}
We note that if $\{\pi_X^i \ \vert \ 0 \leq i \leq 2d\}$ is a generically defined multiplicative Chow--K\"unneth decomposition, then $\{\pi_{X^n}^i \ \vert \ 0 \leq i \leq 2nd\}$ with $\pi_{X^n}^i$ as in \eqref{eq:CK} is a generically defined multiplicative Chow--K\"unneth decomposition for $X^n$; in particular,  $\CH^*(X^n)_{(0)}$ does define a $\QQ$-subalgebra of $ \CH^*(X^n)$.
Finally, the generalized Franchetta conjecture for $X^{n}$
	shows that generically defined cycles on $X^n$ belong to $\CH^*(X^n)_{(0)}$;
	indeed, 
	$(\pi^j_{X^n})_\ast \GDCH^i(X^n)$ consists then of generically defined cycles on $X^n$ which are homologically trivial for $j\neq 2i$.
\end{proof}

\begin{remark}[On assumptions $(a)$, $(b)$ and $(b')$]
	    We note that in Proposition~\ref{prop:F-0} and in Corollary~\ref{cor:F-BVF} below,
		assumption $(b')$ obviously implies assumption $(a)$ and that, provided $n\geq 3$, Proposition~\ref{prop:F-MCK} shows that the combination of $(a)$ and $(b)$ implies assumption $(b')$. In fact, assuming $(a)$ and $(b)$, if $n\geq 3$, then $X$ has a \emph{unique} generically defined Chow--K\"unneth decomposition. 
\end{remark}

\begin{corollary}\label{cor:F-BVF}
Let $\mathcal F$ be the moduli stack of a locally complete family of polarized hyper-K\"ahler varieties and let $X$ be any member of the universal family $\mathcal X \to \mathcal F$. Let $n$ be a positive integer. Assume the following:
\begin{enumerate}[(a)]
\item The very general member of $\mathcal X \to \mathcal F$ satisfies Grothendieck's K\"unneth standard conjecture.
\item 
 The family $\XX^n_{/\mathcal F} \to \mathcal F$ has the Franchetta property.
	\item[(b')]
$X$ has a generically defined multiplicative Chow--K\"unneth decomposition $\{\pi_X^i \ \vert \ 0 \leq i \leq 2d\}$.
\item Murre's conjecture (D) \cite{Murre} holds for the Chow--K\"unneth decomposition on $X^n$ induced as in \eqref{eq:CK} by the generically defined Chow--K\"unneth decomposition on $X$ (which exists and is unique by Proposition~\ref{prop:F-MCK}). 
\end{enumerate}
Then the Beauville--Voisin--Franchetta Conjecture~\ref{conj2} holds for $X^n$.
\end{corollary}
\begin{proof}
	By Proposition~\ref{prop:F-0}, we have the inclusion 
	\[\GDCH^*_B(X^n) \subseteq \CH^*(X^n)_{(0)}\]
	of $\QQ$-subalgebras of $\CH^*(X^n)$, where  $\CH^*(X^n)_{(0)}$ is as in \eqref{eq:CK}.
	Murre's conjecture (D) \cite{Murre} implies that $\pi_{X^n}^{2i}$ acts as zero on the kernel of the cycle class map $\CH^i(X^n) \to H^{2i}(X^n,\QQ)$, in particular that the restriction of the cycle class map to $\CH^i(X^n)_{(0)}$ injects in cohomology.
	To conclude, it remains to see that $\CH^1(X^n)_{(0)} = \CH^1(X^n)$. 
	Since $H^1(X^n,\QQ) = 0$, we have that the cycle class map $\CH^1(X^n) \to H^2(X^n,\QQ)$ is injective. 
	The claim is then a consequence of the fact that the cycle class map is compatible with the action of correspondences; specifically, for our purpose, if  $Y$ is a smooth projective variety such that $\CH^i(Y) \to H^{2i}(Y,\QQ)$ is injective and if $\gamma \in \CH^{\dim Y}(Y\times Y)$ is a correspondence such that $\gamma_*H^{2i}(Y,\QQ) = 0$, then we have $\gamma_*\CH^i(Y)=0$.
\end{proof}

\section{Basic facts in operational Chow cohomology} \label{S:CH}

We collect basic facts related to operational Chow cohomology in the sense of Fulton that are used in the proof of Theorem~\ref{main3}.
For any (possibly singular) variety $M$, we will write $\CH_\ast(M)$ for the Chow groups of~$M$ (with $\QQ$-coefficients), and $\CH^\ast(M)$ for the {\em operational Chow cohomology groups} (with $\QQ$-coefficients) as in \cite[Chapter~17]{F}. By definition, an element in $\CH^i(M)$ is a compatible collection of operations
$\CH_\ast(M^\prime)\to \CH_{\ast-i}(M^\prime)$, for all $M^\prime$ mapping to $M$. The action of an element $a\in\CH^\ast(M)$ on $b\in\CH_\ast(M^\prime)$ is denoted by $a\cap b\in \CH_\ast(M^\prime)$. 
There is a well-defined graded product \cite[17.2($P_1$)]{F} on 
$\CH^\ast(M)$ which is denoted by $a\cdot a^\prime\in\CH^\ast(M)$, and given a morphism of varieties $f\colon N \to M$ there is a well-defined pullback \cite[17.2($P_3$)]{F} denoted by $f^*\colon \CH^\ast(M) \to \CH^\ast(N)$ that commutes with the product \cite[17.2($A_{12}$)]{F}.
If $M$ is of pure dimension~$n$, there is a canonical homomorphism $\CH^i(M) \to \CH_{n-i}(M), \alpha \mapsto \alpha \cap [M]$, which, if $M$ is non-singular, provides a canonical isomorphism $\CH^i(M)\cong\CH_{n-i}(M)$.
The operational Chow cohomology has all the expected formal properties. 
For instance, 
given a fiber square
  \[ \xymatrix{ N^\prime \ar[r]^{g} \ar[d] & M^\prime \ar[d]\\
                           N \ar[r]^{f} & M\\} \]
with $f$ a regular embedding of codimension $e$, and $a\in \CH^\ast(M^\prime)$, one has
  \begin{equation}\label{reg} 
  g^\ast(a)\cap g^\ast(b)= g^\ast(a\cap b)\ \ \ \hbox{in}\ \CH_\ast(N^\prime)\ \ \ \forall\ b\in \CH_\ast(M^\prime)\ ,
  \end{equation}
where $g^\ast\colon \CH_i(M^\prime)\to \CH_{i-e}(N^\prime)$ is the refined Gysin homomorphism of \cite[Chapter~6]{F} (this follows formally from the definition of $\CH^\ast$; see \cite[Definition 17.1($C_3$)]{F}).
Also,
 there is a projection formula: 
 for any proper morphism  $f\colon N \to M$, one has	
 \begin{equation}\label{proj-form}  
 a\cap f_\ast(b)= f_\ast(f^\ast(a)\cap b)\ \ \ \hbox{in}\ \CH_\ast(M)\ ,
 \end{equation}
where $a\in\CH^\ast(M)$ and $b\in\CH_\ast(N)$ \cite[17.2($A_{123}$)]{F};
 in particular, 
if $M$ is smooth and if $f$ is generically finite of degree $d$, then 
the composition 
\begin{equation} \label{comp}
 \CH^i(M)\ \xrightarrow{f^\ast}\ \CH^i(N)\ \xrightarrow{-\cap [N]}\ \CH_{n-i}(N)\ \xrightarrow{f_\ast}\ \CH_{n-i}(M)\cong \CH^i(M)  
\end{equation}
is multiplication by $d$.

	Regarding operational Chow cohomology in degree 1,  there is a natural map (with integral coefficients)
	\[ \pic(M)\ \to\ \CH^1(M)_{\Z}\ .\]
	Indeed, a Cartier divisor on $M$ acts on the Chow groups of any $M^\prime$ mapping to $M$ (by pulling back and intersecting), hence defines an element in the operational Chow cohomology group of codimension 1 with integral coefficients.
	This map is in general not an isomorphism \cite[Example~17.4.9]{F}.

Finally, we will need the following birational invariance property due to Kimura:
\begin{proposition}[Kimura, {\cite[Proposition~3.11]{Kim}}]\label{0cycles}
	The group $\CH^n(-)$ is a birational invariant among (possibly singular) proper $n$-dimensional varieties. More precisely, if $\phi\colon N\dashrightarrow M$ is a birational map among proper $n$-dimensional varieties, there is a canonical isomorphism
	\[ \phi^\ast\colon\ \ \CH^n(M)\ \xrightarrow{\cong}\ \CH^n(N).\]
\end{proposition}

\section{Abel--Jacobi trivial cycles}\label{S:AJ}

The aim of this section is to introduce Abel--Jacobi trivial cycle classes for (not necessarily smooth) projective varieties and show that in codimension 2 these provide a birational invariant (Proposition~\ref{birat}). This will be used in the proof of Theorem~\ref{main3} in the case of codimension-2 cycles. 

The notion of Abel--Jacobi trivial cycle classes on smooth projective varieties is classical, while in the singular case we propose the following (admittedly ad hoc) definition:

\begin{definition}[Abel--Jacobi trivial cycle classes] \label{AJ}  \hfill

(\rom1)
For a smooth projective variety $X$, we define $\CH^\ast_{AJ}(X)$ as the kernel of the Abel--Jacobi map from the subgroup consisting of homologically trivial cycle classes to the Griffiths intermediate jacobian.

(\rom2) For a (possibly singular) projective variety $X$, we define
  \[  \CH_{\ast}^{AJ}(X):= \ima\bigl(  \CH_{\ast}^{AJ}(\wt{X})\to \CH_{\ast}(X)\bigr)\ \ \ \subset\ \CH_{\ast}(X)\ ,\]
  where $\wt{X}\to X$ is a resolution of singularities, i.e., a proper surjective morphism with $\wt{X}$ smooth.
\end{definition}

 Definition~\ref{AJ}(ii) is unambiguous:
\begin{proposition}
Let $X$ be a projective variety. Then $\CH_{\ast}^{AJ}(X)$ does not depend on the choice of a resolution $\wt{X} \to X$.
\end{proposition}
\begin{proof}
Since two given resolutions  of $X$ can be dominated by a third one, 
it is enough to show that given any dominant morphism $f\colon V \to W$ of smooth projective varieties, 
the pushforward $f_\ast\colon \CH_{\ast}(V) \to \CH_{\ast}(W)$ sends $\CH_{\ast}^{AJ}(V)$ surjectively onto $\CH_{\ast}^{AJ}(W)$. The latter follows from the fact that 
$f_*f^* = \deg(f)\Delta_V$ as correspondences (by the projection formula),
and from the fact \cite[Example 1.11]{Saito} 
 that Abel--Jacobi equivalence defines an adequate equivalence relation (in the sense of Samuel); in particular, that correspondences between smooth projective varieties respect the subgroup $\CH_\ast^{AJ}\subset \CH_\ast$. 
\end{proof}

\begin{remark}[Functoriality properties of $\CH_{\ast}^{AJ}$]\label{functAJ} 
We will need the following two functoriality properties of $\CH_{\ast}^{AJ}$.
Let $f\colon X \to Y$ be a morphism of projective varieties.

(\rom1) As is clear from the definition, the proper pushfoward map $f_\ast \colon \CH_{\ast}(X) \to \CH_{\ast}(Y)$ induces a pushforward $f_\ast \colon \CH_{\ast}^{AJ}(X) \to \CH_{\ast}^{AJ}(Y)$.

(\rom2) If $X$ is equidimensional and if $Y$ is smooth, the pullback map
  \[\CH^\ast(Y)  \ \xrightarrow{f^\ast}\ \CH^\ast(X)\ \xrightarrow{\cap [X]}\ \CH_{\dim X-\ast}(X) \]
   sends $\CH^\ast_{AJ}(Y)$ to $\CH^{AJ}_{\dim X-\ast}(X)$. 
   Indeed, if $\tau\colon\wt{X}\to X$ is a resolution of singularities, then the projection formula implies that the indicated pullback map is the same map as
  \[  \CH^\ast(Y) \ \xrightarrow{(f\circ\tau)^\ast}\ \CH^\ast(\wt{X})=\CH_{\dim X-\ast}(\wt{X}) \ \xrightarrow{\tau_\ast} \CH_{\dim X-\ast}(X)\ . \]
But the pullback along $f\circ\tau$  preserves $\CH^\ast_{AJ}$ since Abel--Jacobi equivalence is an adequate equivalence relation.
   \end{remark}

\begin{proposition}\label{birat}
Let $\tau\colon \wt{X}\to X$ be a birational morphism between projective varieties, where $X$ is smooth and equidimensional of dimension $n$. 
The pushforward along $\tau$ induces an isomorphism
  \[ \tau_\ast\colon\ \ \CH_{n-2}^{AJ}(\wt{X})\ \xrightarrow{\cong}\ \CH_{n-2}^{AJ}(X)\ .\]
\end{proposition}

\begin{proof} Let $\wt{\tau}\colon\wt{\wt{X}}\to \wt{X}$ be a resolution of singularities. The composition
  \[ \CH_{n-2}^{AJ}(\wt{\wt{X}})\ \xrightarrow{\wt{\tau}_\ast}\   \CH_{n-2}^{AJ}({\wt{X}})\ \xrightarrow{\tau_\ast}\  \CH_{n-2}^{AJ}(X)  \]
  is an isomorphism (birational invariance between smooth projective varieties; see e.g.\ \cite[Proposition~5.4]{Via-ACF}). This immediately implies the proposition.
  \end{proof}

\section{LLSS eightfolds and the resolution of Voisin's rational map}\label{S:LLSS}

We start by briefly recalling the geometric construction of Lehn--Lehn--Sorger--van Straten~\cite{LLSS}. 
Let $Y \subset \PP^5$ be a smooth cubic fourfold not containing a plane. Twisted cubic curves on $Y$ belong to an irreducible component $M_3(Y)$ of the Hilbert scheme $\textrm{Hilb}^{3m+1}(Y)$; in particular, $M_3(Y)$ is a ten-dimensional smooth projective variety and it is referred to as the Hilbert scheme of \emph{generalized twisted cubics} on $Y$. There exists a hyper-K\"ahler eightfold $Z = Z(Y)$ and a morphism $u: M_3(Y) \to Z$ which factorizes as $u = \Phi \circ a$, where $a: M_3(Y) \to Z'$ is a $\PP^2$-bundle and $\Phi: Z' \to Z$ is the blow-up of the image of a Lagrangian embedding $j: Y \hookrightarrow Z$. 
The hyper-K\"ahler eightfold $Z$ 
will be called an \emph{LLSS eightfold}; by \cite{AL}, it 
is of $K3^{[4]}$-type.

Given a smooth cubic fourfold $Y \subset \PP^5$ not containing a plane, 
the LLSS eightfold $Z = Z(Y)$ is closely related to the Fano variety $F=F(Y)$  of lines in $Y$. Indeed, Voisin~\cite[Proposition~4.8]{V14} constructed a degree $6$ dominant rational map
\begin{equation}\label{psi}
\psi\colon\ \ F\times F\ \dashrightarrow\ Z\ .
\end{equation} 
For the sake of completeness, we briefly recall the geometric construction of the rational map $\psi$. 
Let $(l,l') \in F \times F$ be a generic point, so that the two lines $l, l'$ on $Y$ span a three-dimensional linear space in $\mathbb{P}^5$. For any point $x \in l$, the plane $\langle x, l'\rangle$ intersects the smooth cubic surface $S = \langle l, l' \rangle \cap Y$ along the union of $l'$ and a residual conic $Q'_x$, which passes through $x$. Then, $\psi(l,l') \in Z$ is the point corresponding to the two-dimensional linear system of twisted cubics on $S$ linearly equivalent to the rational curve $l \cup_x Q'_x$ (this linear system actually contains the $\mathbb{P}^1$ of curves $\left\{ l \cup_x Q'_x \mid x \in l\right\}$).
\medskip

For a polarized hyper-K\"ahler variety $(M,h)$, we define $\CH^1(M)_{\rm prim}$ to be the orthogonal complement (with respect to the polarization form or the Beauville--Bogomolov form -- they are the same) of the class of $h$ in $\CH^1(M)$.
Using Voisin's rational map $\psi$, one can already establish the following weak version of Beauville's conjecture :

\begin{proposition}\label{Beauville-weak}
	 Let $Z$ be an LLSS eightfold. The $\QQ$-subalgebra
	\[  \langle \CH^1(Z)_{\rm prim} \rangle\ \ \subset\ \CH^\ast(Z) \]
	generated by the primitive divisors injects into cohomology.
\end{proposition}

\begin{proof} 
	Recall from \cite{Lehn} that there exists an anti-symplectic involution $\iota$ on every LLSS eightfold $Z$; see also \cite[Section 1.1]{FM+}.
	Moreover, as explained in \cite[Section~2.6]{CCL}, the Voisin map
	$\psi\colon F\times F  \dashrightarrow Z$
	is $G\coloneqq \Z/2\Z$-equivariant, where the action of $G$ on~$Z$ is via the 
	involution $\iota$, and the action of $G$ on $X\coloneqq F\times F$ is by switching the two copies of $F$.	
	From equivariant resolution of indeterminacy \cite[Theorem 1]{RY}, 
    there exists a $G$-equivariant morphism $\tau\colon \wt{X}\to X$ which is a composition of blow-ups with smooth $G$-invariant centers, such that the composed map $\wt{\psi}:=\psi\circ\tau$ is a $G$-equivariant morphism.    
	By the blow-up formula for Chow groups, we have
	\[  \CH^1(\wt{X})= \tau^\ast \CH^1(X) \oplus \QQ[E_1]\oplus \cdots \oplus \QQ[E_r] \ ,\]
	where the $E_j$ are the exceptional divisors.
	Since each $E_j$ is $G$-invariant, we obtain
	\[  \CH^1(\wt{X})^+ = \tau^\ast \CH^1(X)^+ \oplus \QQ[E_1]\oplus \cdots \oplus \QQ[E_r]\ ,\ \ \CH^1(\wt{X})^- =   \tau^\ast \CH^1(X)^-\ ,\]
	where the superscript $+$ (resp.\ $-$) indicates the invariant (resp.\ anti-invariant) part under the action of the involution.
	As for $Z$, the decomposition in eigenspaces of $\CH^1(Z)$ is such that
\[ \CH^1(Z)^+ = \QQ[h]\ ,\ \ \CH^1(Z)^- =\CH^1(Z)_{\rm prim}\ ,\]
where $h$ is the polarization on $Z$ coming from that on the cubic fourfold $Y$.
Indeed, $\CH^1(Z)$ injects equivariantly in $H^2(Z,\QQ)$, and $H^2(Z,\QQ)^- = H^2(Z,\QQ)_{\rm prim}$ for a very general member of the family and so for any member of the family.
Hence,
\[    {\wt{\psi}}^\ast \CH^1(Z)_{\rm prim} = {\wt{\psi}}^\ast \CH^1(Z)^- \ \subset\ \CH^1(\wt{X})^- =   \tau^\ast \CH^1(X)^-\ .\]	       	  
We are thus reduced to showing that the $\QQ$-algebra $\langle \CH^1(X)^- \rangle$ injects into cohomology, which follows from Voisin's result \cite{V17} stating that the $\QQ$-algebra $\langle \CH^1(F) \rangle$ injects into cohomology.
\end{proof}

Moving forward, we will need the following explicit resolution of Voisin's map $\psi$.
Its indeterminacy locus was showed by Muratore \cite[Theorem 1.2]{Mura}  to coincide with the locus $I\subset F\times F$ of intersecting lines\,; more precisely, we have:

\begin{proposition}[Chen \cite{Chen}]\label{chen} 
Let $Y$ be a smooth cubic fourfold not containing a plane, and let $\psi$ be the rational map as above. The indeterminacy of $\psi$ is resolved by blowing up  the codimension-$2$ locus $I\subset F\times F$ of intersecting lines. 
That is, let $\wt{X}$ be the blow-up of $X\coloneqq F\times F$ with center $I \coloneqq \{(l,l')\in F \times F \ \vert \ l\cap l' \neq \varnothing\}$. Then there is a morphism $\wt{\psi}$ fitting into a commutative diagram
  \begin{equation*}
\xymatrix{
\wt{X} \ar[d]_{\tau} \ar[dr]^{\wt{\psi}} & \\
X\coloneqq F^{}\times_{}F^{}  \ar@{-->}[r]_-{\psi}& Z.
}
\end{equation*}  
\end{proposition}

\begin{proof} This is \cite[Theorem 4.6(a)]{Chen}, which uses the interpretation of $Z$ and $F$ as moduli spaces of stable objects in the Kuznetsov component of $Y$.
\end{proof}

\begin{remark}
	One could have used the resolution of $\psi$ provided by Proposition~\ref{chen} as a replacement to the abstract equivariant resolution of indeterminacy of \cite{RY} to prove Proposition~\ref{Beauville-weak}.
	This however, for cubic fourfolds for which the blow-up $\wt{X}$ of Proposition~\ref{chen} is singular,
	would have come at the cost of using the cycle class map of Bloch--Gillet--Soul\'e \cite{BGS} from operational Chow cohomology
	$ \CH^i(\wt X)$ to  $\Gr_{2i}^W H^{2i}(\wt X,\QQ)$, where $W_\ast$ refers to Deligne's weight filtration.
\end{remark}

\begin{lemma}\label{jiang} Let $Y$ be a smooth cubic fourfold not containing a plane. 
	Let $X\coloneqq F(Y)\times F(Y)$ and $\tau\colon \wt{X}\to X$ be as in Proposition \ref{chen}, and let $E\subset\wt{X}$ denote the exceptional divisor of $\tau$. Then
  \[ \pic(\wt{X})_{\QQ} =  \tau^\ast \pic(X)_{\QQ}\oplus \QQ[E]\ .\]
\end{lemma}

\begin{proof} 
Let $I\subset X=F\times F$ denote the locus of intersecting lines, as in Proposition~\ref{chen}. 
We observe that both the blow-up $\wt{X}$ and the exceptional divisor $E$ of the blow-up are irreducible. Indeed, Jiang \cite[Section~4.3.3]{Jiang} has computed the codimensions of the degeneracy loci $X^{\ge i}({\mathcal I}_I)$ (in the notation of loc.\ cit.): the codimension is 6 for $i=3$ and the degeneracy locus is empty for $i\ge 4$. The irreducibility of $\wt{X}$ and $E$ then follow from a general result of Ellingsrud--Str\o mme \cite[Proposition~3.2]{ES}.

The blow-up of Proposition \ref{chen} gives rise to an exact sequence of Chow groups
$$\xymatrix{\CH_i(E)\ar[r] & \CH_i(\wt{X})\oplus \CH_i(I)\ar[r] & \CH_i(X)\ar[r]& 0\ ,
}$$
  where all arrows are induced by pushforward morphisms.
  Such a sequence is obtained by a diagram chase on  the following commutative diagram with exact rows
  $$\xymatrix{\CH_i(U,1) \ar@{=}[d] \ar[r] & \CH_i(E) \ar[d]^{\tau_{E,\ast}}\ar[r] & \CH_i(\wt{X}) \ar[r] \ar[d]^{\tau_\ast} & \CH_i(U) \ar@{=}[d] \ar[r] & 0 \\
  	\CH_i(U,1)  \ar[r] & \CH_i(I) \ar[r] & \CH_i({X}) \ar[r]  & \CH_i(U) \ar[r] & 0   	
  }$$
  where $U = \wt{X} \setminus E = X \setminus I$, and $\CH_i(-,1)$ denotes Bloch's higher Chow groups with $\QQ$-coefficients \cite{Bl1, Bl2}. (Alternately, one can use Voevodsky's bigraded motivic cohomology.)
  
  In particular, taking $i=7$ (and
  using that $E$ is irreducible), we find that the kernel of
  \[ \tau_\ast\colon\ \CH_7(\wt{X})\ \to\ \CH_7(X) \]
  is one-dimensional, generated by the class of $E$.
  Since $X$ is smooth, the composition $\tau_\ast\tau^\ast$ in the commutative diagram
  $$\xymatrix{\pic(X) \ar[r]^{\tau^\ast} \ar@{=}[d] & \pic(\wt{X}) \ar[d] \\
  	\CH_7(X)_{\Z} \ar[r]^{\tau^\ast} & \CH_7(\wt{X})_{\Z} \ar[r]^{\tau_\ast} & \CH_7(X)_{\Z}
  }$$
   is an isomorphism. Here, $\CH_7(-)_\Z$ denotes the Chow group of 7-dimensional cycles with integral coefficients.
   
 Now, we claim that the variety $\wt{X}$ is normal.
 Observe indeed that $\wt{X}$ is smooth outside of  $\tau^{-1}(\hbox{Sing} (I))$, which has codimension at least 2. Since the blow-up morphism $\tau$
 is lci 
 (by virtue of the fact that the ideal of $I$ has homological dimension 1, cf.~\cite[Proof of Proposition~4.5]{Chen}),  $\wt{X}$ is also lci, hence Cohen--Macaulay. The claim thus follows from Serre's criterion for normality.
 
 The normality of $\wt{X}$ implies that 
  the natural map 
   \[ \pic(\wt{X})\ \to\ \CH_7(\wt{X})_{\Z} \]
   is injective. It follows that
   \[ \tau^\ast\colon\ \ \pic(X)\ \to\ \pic(\wt{X}) \]
   is injective with cokernel of rank at most 1. Since the exceptional divisor $E$ (which is an effective Cartier divisor, e.g., by \cite[\href{https://stacks.math.columbia.edu/tag/02OS}{Tag 02OS}]{stacks-project}) is not in the image of $\tau^\ast$, this proves the corollary.
   \end{proof}

   Note that for {\em general\/} cubics $Y$, the blow-up $\wt{X}$ is smooth (see  \cite[Section~5.3]{Jiang}) and in this case the relation of Lemma \ref{jiang} is classical. 
   The remarkable thing
   is that the relation of Lemma \ref{jiang} continues to hold even when $\wt{X}$ is singular.

\section{Proof of Theorem~\ref{main3}}\label{S:main}

Let $B \subset \PP (H^0(\PP^5, \mathcal{O}(3)))$ be the Zariski-open subset parameterizing smooth complex cubic fourfolds not containing a plane. Let $\YY\to B$, $\FFF\to B$ and $\ZZZ\to B$
denote the corresponding universal families of cubic fourfolds, resp. their associated Fano varieties of lines, resp.\ the associated LLSS eightfolds.
Recall from Definition~\ref{def:gdc} that $\GDCH^i_B(-)$ denotes the generically defined cycles.

\begin{theorem}\label{main} Let $Y\subset\PP^5(\C)$ be a smooth cubic fourfold not containing a plane, and let $Z\coloneqq Z(Y)$ be the associated Lehn--Lehn--Sorger-van Straten eightfold.
Let
   \[ R_B^\ast(Z)\coloneqq\Bigl\langle  \GDCH^\ast_B(Z), \ \CH^1(Z)\Bigr\rangle\ \ \subset\ \CH^\ast(Z) \]
   be the $\QQ$-subalgebra generated by generically defined cycles and divisors. The cycle class map induces injections
     \[ R^i_B(Z)\ \hookrightarrow\ H^{2i}(Z,\QQ)\ \ \ \hbox{for}\ i=2\ \hbox{and\ for}\ i=8\ .\]
\end{theorem}

\begin{proof}[Proof of Theorem~\ref{main} in case $i=8$.]
The diagram of Proposition \ref{chen} exists in a family-version, i.e., there is a diagram of schemes over $B$
\begin{equation*}
\xymatrix{
\wt{\XX} \ar[d]_{\tau_{\XX}} \ar[dr]^{\wt{\psi}} & \\
\XX\coloneqq\FF^{}\times_{B^{}}\FF^{}  \ar@{-->}[r]_-{\psi}& \ZZZ.
}
\end{equation*}
In particular, it makes sense to talk about generically defined cycles (with respect to $B$) on $\wt{X}$ and on $X$. We note that $\wt{X}$ may in general be singular (it is known to be smooth only for sufficiently general $Y$, cf. \cite[Section~5.3]{Jiang}), and so some care is needed. We define
  \[ \GDCH^\ast_B(\wt{X})\coloneqq\operatorname{im}\Bigl(  \CH^\ast(\wt{\XX})\to \CH^\ast(\wt{X})\Bigr)\ ,\]
  where $\CH^\ast$ denotes Fulton--MacPherson's operational Chow cohomology (cf. \cite{F} and Section \ref{S:CH} above).

\vskip0.2cm
{\it First step: reduction to $\wt{X}$.} Let us define the $\QQ$-subalgebra
  \[ R_B^\ast(\wt{X})\coloneqq \Bigl\langle   \GDCH^\ast_B(\wt{X}), \ \operatorname{im}\bigl( \pic(\wt{X})\to \CH^1(\wt{X})\bigr)\Bigr\rangle\ \ \subset\ \CH^\ast(\wt{X}) \ .\]
We claim that if $R_B^8(\wt{X})$
 is one-dimensional, then so is  $R_B^8(Z)$.
Indeed, it is readily seen that the pullback map $\wt{\psi}^\ast\colon \CH^\ast(X)\to \CH^\ast(\wt{X})$ sends $R_B^\ast(X)$ to $R_B^\ast(\wt{X})$ (since $\wt{\psi}^\ast$ sends generically defined cycles to generically defined cycles, Cartier divisors to Cartier divisors, and intersection products to intersection products). The claim then follows from the fact that the composition
  \[ \CH^8(Z)\ \xrightarrow{\wt{\psi}^\ast}\ \CH^8(\wt{X})\ \xrightarrow{ - \cap[\wt{X}]}\ \CH_0(\wt{X})\     \xrightarrow{\wt{\psi}_\ast}\ \CH_0(Z)=\CH^8(Z) \]
 is multiplication by $\deg \wt{\psi} = \deg \psi = 6$; see \eqref{comp}.

\vskip0.2cm
{\it Second step: reduction to ${X}$.} Let us consider the $\QQ$-subalgebra
  \[ R_B^\ast({X})\coloneqq \Bigl\langle   \GDCH^\ast_B({X}), \  \CH^1({X})\Bigr\rangle\ \ \subset\ \CH^\ast({X}) \ .\]  
  We claim that $R_B^8(\wt{X})$
is one-dimensional, provided $R_B^8(X)$ is one-dimensional.
  Indeed, we know from Lemma~\ref{jiang} and its proof that
    \[  \begin{split} \operatorname{im}\bigl( \pic(\wt{X})_\QQ \to \CH^1(\wt{X})\bigr)  &= \operatorname{im}\bigl(  \CH^1(X) \xrightarrow{\tau^\ast} \CH^1(\wt{X})\bigr) \oplus \QQ[E]  \\     
    &\subseteq \operatorname{im}\bigl(  \CH^1(X) \xrightarrow{\tau^\ast} \CH^1(\wt{X})\bigr) + \GDCH^1_B(\wt{X}) \ .\\
    \end{split} \]
    As such, $R_B^\ast(\wt{X})$ can be rewritten as
    \[ R_B^\ast(\wt{X})=\Bigl\langle  \GDCH^\ast_B(\wt{X}), \ \operatorname{im}\bigl(  \CH^1(X) \xrightarrow{\tau^\ast} \CH^1(\wt{X})\bigr) \Bigr\rangle\ .\]
    That is, any element $\alpha\in R_B^\ast(\wt{X})$ can be written as a sum of elements of the form
    \[  \alpha= \beta\cdot \tau^\ast( D_1\cdots D_r)\ \ \hbox{in}\ \CH^\ast(\wt{X})\ ,\]
    where $\beta\in\GDCH^\ast_B(\wt{X})$ is generically defined, and $D_1,\ldots, D_r$ are divisors on $X$. 
    Taking the push-forward to $X$, we find 
      \[ \begin{split} \tau_\ast(\alpha\cap [\wt{X}])&= \tau_\ast( \beta\cdot \tau^\ast( D_1\cdots D_r)\cap [\wt{X}])\\
                                                                          &= \tau_\ast \bigl( \tau^\ast(D_1\cdots D_r)\cap (\beta\cap[\wt{X}])\bigr)\\
                                                                          &=(D_1\cdots D_r)\cap \tau_\ast(\beta\cap [\wt{X}]) \\
                                                                          &=(D_1\cdots D_r)\cdot \tau_\ast(\beta\cap [\wt{X}])\ \ \ \hbox{in}\ \CH_\ast(X)\ .
      \end{split}\]
       Here the second equality expresses that operational Chow cohomology acts on Chow groups, the third equality follows from \eqref{proj-form},
       and the fourth is a change of notation reflecting that $X$ is smooth. 
       Since $\beta\in\GDCH^\ast_B(\wt{X})$ is generically defined,
         the element $\tau_\ast(\beta\cap [\wt{X}])\in\CH_\ast(X)$ is also 
       generically defined: indeed, this follows from the diagram
       \[ \xymatrix{  \CH^\ast(\wt{\XX}) \ar[r] \ar[d]^{\cap [\wt{\XX}]} &\CH^\ast(\wt{X})  \ar[d]^{\cap [\wt{X}]}   \\
                                \CH_\ast(\wt{\XX}) \ar[r] \ar[d]^{(\tau_{{\XX}})_\ast} &\CH_\ast(\wt{X})  \ar[d]^{\tau_\ast}  \\
 \CH_\ast({\XX}) \ar[r]  &\CH_\ast({X})   \\  } \]
 where the upper square is commutative by \eqref{reg}, and the lower square is commutative because refined Gysin homomorphisms commute with proper pushforward \cite[Theorem 6.2]{F}. 
 We may thus conclude that      
    \[  \tau_\ast(\alpha\cap [\wt{X}])\ \ \in\ R_B^\ast(X)\ .\]
    Since the composite map
      \[ \CH^8(X) \ \xrightarrow{\tau^\ast}\ \CH^8(\wt{X}) \ \xrightarrow{  - \cap[\wt{X}]}\ \CH_0(\wt{X})\     \xrightarrow{\tau_\ast}\ \CH_0(X)=\CH^8(X)\]
     is the identity by the projection formula (equality \eqref{proj-form}), and since $\tau^*$ is an isomorphism (Proposition~\ref{0cycles}), this proves the claim of the second step.                                                                                                   
                                                                                                        
\vskip0.2cm
{\it Third step: proving the result for ${X} = F\times F$.} 
To prove the theorem, it now remains to prove that $R_B^8(F\times F)$ is one-dimensional, where $R_B^\ast(F\times F)$ is as defined in the second step.

Let us first consider the subring
$$R_B^\ast(F) \coloneqq \Bigl\langle   \GDCH^\ast_B({F}), \  \CH^1({F})\Bigr\rangle\ \ \subset\ \CH^\ast(F) \ .$$
The Fano variety of lines $F$ is contained in the Grassmannian $\mathrm{Gr}(\PP^1,\PP^5)$ parameterizing lines in $\PP^5$ and we denote by $g \in \CH^1(F)$ and $c\in \CH^2(F)$ the first and second Chern classes, respectively, of the rank-2 quotient bundle induced on $F$.
By \cite[\S 3]{FLV} together with the identity $c_2(F) = 5g^2-8c$ \cite[Lemma~A.1]{SV}, we have $$\GDCH^\ast_B({F}) = \langle g, c \rangle.$$
Hence 
$$R_B^\ast(F) = \bigl\langle   c, \  \CH^1({F})\bigr\rangle\ \ \subset\ \CH^\ast({F}) \ .$$
By Voisin~\cite[Theorem~1.4$(ii)$]{V17}, $R_B^\ast(F)$ injects in cohomology via the cycle class map.
Moreover, we have the following relations: 
by \cite[Lemma~3.5 and (3.36)]{V17} we have $c \cdot \CH^1(F) = \QQ[g^3]$,  by \cite[Lemma~3.9]{V17} we have $\CH^1(F)^{\cdot 3} = g^2 \cdot \CH^1(F)$, and by \cite[p.19]{V17} there is a point $o_F$ such that $R_B^4(F) = \QQ[o_F]$.

On the other hand, by \cite[Proposition~6.3]{FLV}, we have 
$$\GDCH^\ast_B({F}\times F) = \langle \Delta, I, g_1, c_1, g_2, c_2 \rangle\ ,$$
where $g_i = p_i^*g, c_i = p_i^*c$ for $p_i : F \times F \to F$ the projection on the $i$-th factor, $\Delta$ is the class of the diagonal inside $F\times F$ and $I$ is the \emph{incidence correspondence}, i.e., the class of the codimension-2 locus of intersecting lines as in Proposition~\ref{chen}. 
Hence
\begin{equation}\label{RFF} R_B^\ast(F\times F) = \Bigl\langle   \Delta, I, c_1, c_2, D_1, D'_2 \ \vert \ D, D' \in \CH^1({F})\Bigr\rangle\ \ \subset\ \CH^\ast({X})\ ,\end{equation}
where again $D_1 \coloneqq p_1^\ast D$ and $D'_2 \coloneqq p_2^\ast D'$.
In order to conclude, it suffices to show that any (weighted) monomial of degree 8 in the variables $\Delta_F, I, c_1, c_2, D_1, D'_2$ with $D, D'$ running through codimension-1 cycles in $\CH^1({F})$ is a multiple of the class of $o_F \times o_F$. 
Let us denote by $(\mathrm{dec})$ any decomposable cycles, i.e., any cycle in the subring 
$$R_{\mathrm{dec}}^\ast(F\times F) \coloneqq R_B^\ast(F)\otimes R_B^\ast(F) = \Bigl\langle    c_1, c_2, D_1, D'_2 \ \vert \ D, D' \in \CH^1({F})\Bigr\rangle\ \ \subset\ \CH^\ast({X}),$$
We have the following relations in $R_B^\ast(F\times F)$:
\begin{enumerate}[(i)]
	\item $I\cdot \Delta = \Delta\cdot (\mathrm{dec})$ by  \cite[Proposition~17.5]{SV}.
	\item $I^2 = 2 \Delta + I\cdot  (\mathrm{dec}) +  (\mathrm{dec})$ by \cite{V17} and \cite[Proposition~17.4]{SV}.
\end{enumerate}
From these two relations, it follows that any monomial of degree 8 is a $\QQ$-linear combination of monomials of the form 
$$(\mathrm{dec}), \ \Delta^2, \ \Delta\cdot (\mathrm{dec}), \ \mbox{and}\ I \cdot (\mathrm{dec}).$$
Moving forward, first note that $R_{\mathrm{dec}}^8(F\times F) = R_B^4(F) \otimes R_B^4(F) = \QQ [o_F \times o_F]$, so that $\mathrm{(dec)}$ lies in $\QQ [o_F \times o_F]$.
Second, 
denoting by $\iota_{\Delta}\colon F \to F\times F$ the diagonal embedding, we have
$\Delta^2 = \iota_{\Delta,\ast}\iota_{\Delta}^\ast \Delta$ and since $\iota_{\Delta}^\ast \Delta$ belongs to $\GDCH^4_B(F) = \QQ [o_F]$, we find that $\Delta^2$ is a multiple of the class of $o_F\times o_F$.
Third, since $\Delta \cdot \alpha_2 = \alpha_1 \cdot \Delta$ for any $\alpha\in \CH^\ast(F)$, we have that $\Delta\cdot (\mathrm{dec})$ belongs to $\Delta \cdot p_2^\ast R_B^4(F) = \Delta \cdot \QQ [p_2^\ast o_F] = \QQ[o_F \times o_F]$.
Finally, from $R_B^4(F) = \QQ[o_F] = \QQ[c^2]$, from $c \cdot \CH^1(F) = \QQ[g^3]$ and from $c_1\cdot I = (\mathrm{dec})$ and $I\cdot c_2 = (\mathrm{dec})$ (see~\cite[Lemma~17.6]{SV}), we find that $I \cdot (\mathrm{dec})$ is a $\QQ$-linear combinations of cycles of the form 
$$ (\mathrm{dec}) \ \  \mbox{and} \ \ \alpha_1 \cdot I \cdot \beta_2, \quad \mbox{for}\ \alpha, \beta \in \CH^1(F)^{\cdot 3}.$$
Using that  $\CH^1(F)^{\cdot 3} = g^2 \cdot \CH^1(F)$, we get that 
$\alpha_1 \cdot I \cdot \beta_2$ is a $\QQ$-linear combination of terms of the form
$$D_1\cdot g_1^2 \cdot I \cdot g_2^2 \cdot D'_2, \quad \mbox{for some } D,D'\in \CH^1(F).$$
By \cite[Proposition~6.4]{FLV}, $g_1^2\cdot I \cdot g_2^2$ is a $\QQ$-linear combination of cycles of the form $(\mathrm{dec})$ and $g_1^2\cdot \Delta$. (The latter is based on the relation \cite[Theorem~A.1]{FLV}.) 
Therefore, we find that $\alpha_1 \cdot I \cdot \beta_2$ is a $\QQ$-linear combination of decomposable cycles, and we conclude using that 
$R_{\mathrm{dec}}^8(F\times F) = R_B^4(F) \otimes R_B^4(F) = \QQ [o_F \times o_F]$.
This proves the $i=8$ part of Theorem~\ref{main}.
\end{proof}

\begin{proof}[Proof of Theorem~\ref{main} in case $i=2$.]
We will adhere to the notations and the structure of the above three-step argument for $i=8$.

In a first step, we reduce from $Z$ to $\wt{X}$. Let us consider 
  \[  R_i(\wt{X}):= \ima\Bigl( R^{8-i}_B(\wt{X})\to \CH_i(\wt{X})\Bigr)\ \ \subset\ \CH_i(\wt{X}) \]
  (note that $R_\ast(\wt{X})$ is not a ring; it is a graded $R^\ast_B(\wt{X})$-module).
We claim that in order to prove the injectivity of $R^2_B(Z)\to H^4(Z,\QQ)$, it suffices to prove the vanishing 
  \[    R_6^{AJ}(\wt{X}):=R_6(\wt{X})\cap \CH^{AJ}_6(\wt{X})=0 \]
   (where $\CH_6^{AJ}$ is as in Definition \ref{AJ}). Indeed, 
the composition
  \[ \CH^2(Z)\ \xrightarrow{\wt{\psi}^\ast}\ \CH^2(\wt{X})\ \xrightarrow{ - \cap[\wt{X}]}\ \CH_6(\wt{X})\     \xrightarrow{\wt{\psi}_\ast}\ \CH_6(Z)=\CH^2(Z) \]
 is multiplication by $\deg \wt{\psi} = \deg \psi = 6$ (see equation \eqref{comp}) and preserves the subgroup $\CH_\ast^{AJ}$ (Remark \ref{functAJ}). Thus, if $R_6^{AJ}(\wt{X})=0$, then 
 also $R^2_B(Z)\cap \CH_6^{AJ}(Z)=0$. 
 But the hyper-K\"ahler variety $Z$ has no odd-degree cohomology, and so $\CH_\ast^{AJ}(Z)=\CH_\ast^{hom}(Z)$, thereby proving the claim.
 
 In a second step, we go from $\wt{X}$ to $X$. We claim that in order to prove the vanishing $R_6^{AJ}(\wt{X})=0$, it suffices to prove that $R^2_B(X)$ injects into cohomology under the cycle class map. Indeed, we have already seen (in the second step for $i=8$ above) that the pushforward $\tau_\ast$ sends $R_\ast(\wt{X})$ to $R^\ast_B(X)$. In addition, we know that  the composite map
      \[ \CH^2(X) \ \xrightarrow{\tau^\ast}\ \CH^2(\wt{X}) \ \xrightarrow{  - \cap[\wt{X}]}\ \CH_6(\wt{X})\     \xrightarrow{\tau_\ast}\ \CH_6(X)=\CH^2(X)\]
     is the identity by \eqref{comp}. Since the restriction of  $\tau_*$ to Abel--Jacobi trivial cycles is an isomorphism (Proposition~\ref{birat}) and $\CH_\ast^{AJ}(X)=\CH_\ast^{hom}(X)$ (as $X$ has no odd-degree cohomology), this proves the claim of the second step.  
     
In a third and final step, we prove the required result on $X$, i.e. we prove that $R^2_B(X)\cap \CH_6^{hom}(X)=0$. It follows from \eqref{RFF} that
  \[ R^2_B(X)=R^2_{dec}(F\times F) +\QQ[I]\ .\]
  But the class of $I$ is linearly independent of the decomposable classes in cohomology, since, by \cite[Proposition~6]{BD}, it
  induces an isomorphism $I_\ast\colon H^{4,2}(F) \xrightarrow{\cong} H^{2,0}(F)$.                                                                                             
   \end{proof}

\begin{proof}[Proof of Theorem~\ref{main3}]
	
	Theorem~\ref{main} implies Theorem~\ref{main3} stated in the introduction.
	 Indeed, the canonical morphism from the parameter space $B$ to the moduli stack of LLSS eightfolds~$\mathcal F$ is surjective, so that there are inclusions
	\[\GDCH^\ast_{\mathcal F}(Z)\ \subset\ \GDCH^\ast_B(Z)\ \subset\ \CH^\ast(Z)\]
	(see, e.g., \cite[Remark~2.6]{FLV}), and hence
	an inclusion $R^\ast(Z)\subset R^\ast_B(Z)$.
	
	As for the ``moreover'' part of Theorem~\ref{main3}, we observe that the proof of 	Theorem~\ref{main} actually shows that
	\[ R^8_B(Z) = \QQ \psi_\ast(o_F\times o_F)\ ,\]
	where $\psi$ is Voisin's rational map and $o_F\in \CH^4(F)$ is the distinguished point. Since $o_F$ is represented by any point on the constant cycle surfaces $V\subset F$ constructed by Voisin \cite[Proposition 4.5]{V14}, the 1-dimensional space $R^8_B(Z)$ is spanned by any point on the 4-dimensional constant cycle subvariety $W\subset Z$ constructed by Voisin \cite[Corollary~4.9]{V14} as the closure of the image $\psi(V\times V^\prime)$. 
	That is, defining $o_Z\in\CH^8(Z)$ as the class of any point on any of the constant cycle subvarieties $W\subset Z$ constructed by Voisin, we have
	\[ R^8(Z)=R^8_B(Z)=\QQ[o_Z]\ .\]
\end{proof}

\begin{remark}
The Lagrangian embedding of the cubic fourfold $Y$ in $Z=Z(Y)$ provides another 4-dimensional constant cycle subvariety. The embedding $Y\subset Z$ being generically defined, Theorem \ref{main} implies that points of $Y$ also represent the class $o_Z\in\CH^8(Z)$ defined above.
\end{remark}

		 \bibliographystyle{alpha}

\end{document}